\input amstex
\documentstyle{amsppt}
\magnification=\magstep1
\hsize=5.2in
\vsize=6.8in

\centerline  {\bf  ERGODIC SUBEQUIVALENCE RELATIONS INDUCED}
\vskip 0.05in
\centerline {\bf   BY A BERNOULLI ACTION}

\vskip 0.2in
\centerline {\rm by}
\vskip 0.05in
\centerline {\rm IONUT CHIFAN and ADRIAN IOANA}

\address Math Dept, UCLA, LA, CA 90095-155505 and IMAR, Bucharest, Romania
\endaddress
\email ichifan\@math.ucla.edu\endemail

\address Math Dept,  Caltech, Pasadena, 91125 and IMAR, Bucharest, Romania
\endaddress
\email aioana\@caltech.edu \endemail  
\topmatter
\abstract
Let $\Gamma$ be a countable group and denote by $\Cal S$  the equivalence relation induced by the Bernoulli action $\Gamma\curvearrowright [0,1]^{\Gamma}$, where $[0,1]^{\Gamma}$ is endowed with the product Lebesgue measure.  We prove that for any subequivalence relation $\Cal R$ of $\Cal S$, there exists a partition $\{X_i\}_{i\geq 0}$ of $[0,1]^{\Gamma}$ with $\Cal R$-invariant measurable sets such that $\Cal R_{|X_0}$ is hyperfinite and $\Cal R_{|X_i}$ is strongly ergodic (hence ergodic), for every $i\geq 1$. 
\endabstract

\endtopmatter

\document

\head  \S {1. Introduction and statement of results.}\endhead
\vskip 0.1in
During the past decade there have been many interesting new directions arising in the field of measurable group theory.
One direction came from the  {\bf deformation/rigidity theory} developed recently by S. Popa  in order to study group actions and von Neumann algebras ([P5]).  
Using this theory, Popa obtained striking rigidity results concerning the equivalence relations and the II$_1$ factors induced by a Bernoulli action ([P1-4]).

To recall these results, let $\Gamma$ be a countable group and 
let $\Gamma\curvearrowright [0,1]^{\Gamma}$ be the Bernoulli action, where $[0,1]^{\Gamma}$ is endowed
 with the product Lebesgue measure. Denote by $\Cal S$ the induced equivalence relation, i.e.
 $(x_{\gamma})\sim (y_{\gamma})$ if there exists $\gamma'$ such that $x_{\gamma}=y_{\gamma'\gamma}$, for 
all $\gamma$ in $\Gamma$. Popa then  proved that if $\Gamma$  satisfies a strong non-amenability condition
 (e.g. if  $\Gamma$ has property (T) or  splits as a non-amenable product of two infinite groups), then
 $\Cal S$ remembers both the group and the action ([P3,4]). Moreover, he showed  that for any non-amenable
 group $\Gamma$, the associated II$_1$ factor $L(\Cal S)$ is prime ([P4]).

The main goal of this paper is to present a new rigidity phenomenon displayed by the equivalence relation $\Cal S$.  
More precisely, we show that for  any countable group $\Gamma$, we have the following structure result for  the subequivalence relations of $\Cal S$:

\proclaim {Theorem 1} Let $\Cal R\subset\Cal S$ be a subequivalence relation. Then there exists a measurable partition $\{X_i\}_{i\geq 0}$ of $[0,1]^{\Gamma}$ with $\Cal R$-invariant sets such that
\vskip 0.02in
$(a)$. $\Cal R_{|X_0}$ is hyperfinite.
\vskip 0.02in
$(b)$. $\Cal R_{|X_i}$ is strongly ergodic (therefore ergodic), for all $i\geq 1$.
\vskip 0.05in
Moreover, the same holds for any quotient equivalence relation $\Cal R'$ of $\Cal R$.
\endproclaim

\vskip 0.05in 

In particular, Theorem 1 shows that for the non-hyperfinite subequivalence relations of $\Cal S$, the 
notions of {\bf ergodicity} and {\bf strong ergodicity} are equivalent. For a more general statement, see Theorem 7.

To prove Theorem 1, we follow an operator algebra approach using Popa's {\bf deformation/rigidity strategy}. 
In this respect, recall first that to every  countable, measure preserving (m.p.), equivalence relation 
$\Cal R$ one can associate a finite von Neumann algebra $L(\Cal R)$ ([FM]).  
Then an inclusion of equivalence relations $\Cal R\subset\Cal S$ gives an inclusion of von Neumann 
algebras $L(\Cal R)\subset L(\Cal S)$. On the other hand, remark that for our particular
 $\Cal S$,   we can view $L(\Cal S)$ as the  wreath product II$_1$ factor
 $(\overline{\bigotimes}_{\Gamma}L^{\infty}[0,1])\rtimes\Gamma$. Altogether, we get
 that $L(\Cal R)\subset (\overline{\bigotimes}_{\Gamma}L^{\infty}[0,1])\rtimes\Gamma$.  Theorem 1 
 is then  a consequence of the following general result on controlling relative commutants in wreath product factors.

\proclaim {Theorem 2} Let $(B,\tau)$ be an amenable, finite von Neumann algebra and let $\Gamma$ be a countable group. On the infinite tensor product $\overline{\bigotimes}_{\Gamma}B$ consider the Bernoulli action of $\Gamma$ and denote by $M$ the crossed product $(\overline{\bigotimes}_{\Gamma}B)\rtimes\Gamma$.
Let $p\in M$ be a projection, let $P\subset pMp$ be a von Neumann subalgebra with no amenable direct summand and denote $Q=P'\cap pMp$. 
\vskip 0.03in
(i) Then there exists a non-zero partial isometry $v\in M$ such that $v^*v\in Q'\cap pMp$ and $vQv^*\subset L\Gamma$.
\vskip 0.03in
(ii) Moreover, if $\Gamma$ is ICC (infinite conjugacy class), then there exists a unitary $u\in M$ such that $u(P\vee Q)u^*\subset L\Gamma$.

\endproclaim   
Theorem 2 has been proved by Popa (Lemma 5.2. in [P4]) under the additional assumption that no corner of $Q$ embeds into $\overline{\bigotimes}_{\Gamma}B$ by using the {\bf malleable deformations} of Bernoulli actions. In fact, our proof of Theorem 2 follows closely the proof of Lemma 5.2. in [P4]. 
 Indeed, just as in [P4], we use a spectral gap argument to show that since $Q$ is the commutant of a non-amenable algebra, then it behaves, in some sense (i.e. with respect to certain deformations of $M$), as a rigid subalgebra of $M$. 

The main difference in our approach is that we use the {\bf (weakly) malleable deformations} of Bernoulli actions considered  in section 2 of [I] rather than the malleable ones.
The benefit of this approach is that we can apply Theorem 3.3. in [I] to get precise information on the position of the "rigid" subalgebra $Q$ in the wreath product factor $M$.

 We point out that in the proof of Theorem 1 we actually only use the following consequence of Theorem 2. We will later see that  the next Corollary is in fact equivalent to Theorem 1 (Proposition 6).
 
\proclaim {Corollary 3} Assume that $(B,\tau)$ is amenable. Let $p\in M$ be a projection and let $Q\subset p(\overline{\bigotimes}_{\Gamma}B)p$ be a diffuse von Neumann algebra. 

\vskip 0.03in	
(i) Then $Q'\cap pMp$ is amenable.
\vskip 0.03in
(ii) Moreover, if  $N\subset pMp$ is a non-amenable subfactor containing $Q$, then $N'\cap Q^{\omega}=\Bbb C1$.
\endproclaim

In the case $\Gamma$ is exact, part (i) of Corollary 3 has been first proved by N. Ozawa  using C$^*$-algebraic methods (Theorem 4.7. in [O2]). Ozawa's result thus provides a different approach to proving Theorem 1, in the exact case. 

\vskip 0.1in
Aside from the introduction, this paper has three more sections.
The  second section  deals with the technical tools that are needed in the proof of Theorem 2. In Section 3, we first prove Theorem 2 and then deduce Corollary 3 and Theorem 1. Lastly, in Section 4, we note a few  applications of our main results.   

\vskip 0.1in
{\bf Acknowledgment.} We are grateful to Professor Sorin Popa for suggesting to us the idea from which this project arose.
The second author is also grateful to Professor Alekos Kechris for useful discussions.
\head \S {2. Technicalities.}\endhead
\vskip 0.1in 
In this Section we discuss some of the main ingredients of the proof of Theorem 2. Since the context of this proof is the same as in Section 2 in [I], we begin by recalling some notations and constructions from there. Let $(B,\tau)$ be a finite von Neumann algebra with a  normal, faithful trace $\tau$.  For a countable set $I$, we  denote by $\overline{\bigotimes}_{I}B$ the tensor product von Neumann algebra $\overline{\otimes}_{i\in I}(B,\tau)_i$.

Let $(\tilde B,\tilde\tau)$ denote the free product von Neumann algebra of $(B,\tau)$ with the group von Neumann algebra $L\Bbb Z$.
On $\overline{\bigotimes}_{\Gamma}\tilde B$ consider the Bernoulli action of $\Gamma$ given by $\tilde\sigma(\gamma)(\otimes_{\gamma'}x_{\gamma'})=\otimes_{\gamma'}x_{\gamma^{-1}\gamma'}.$ Then $\tilde\sigma$ leaves the subalgebra $\overline{\bigotimes}_{\Gamma}B\subset\overline{\bigotimes}_{\Gamma}\tilde B$ invariant and the restriction $\sigma={\tilde\sigma}_{|\overline{\bigotimes}_{\Gamma}B}$ is precisely the Bernoulli action of $\Gamma$ on $\overline{\bigotimes}_{\Gamma}B$.
Denote by $M$ (resp. $\tilde M$)  the crossed product von Neumann algebra $(\overline{\bigotimes}_{\Gamma}B)\rtimes_{\sigma}\Gamma$ (resp. $(\overline{\bigotimes}_{\Gamma}\tilde B)\rtimes_{\tilde\sigma}\Gamma$). Also,  denote by $\{u_{\gamma}\}_{\gamma\in\Gamma}\subset \tilde M$ the canonical unitaries implementing $\tilde\sigma$.

Next, let $u\in L\Bbb Z$ be the generating Haar unitary and let $h=h^*\in L\Bbb Z$ such that $u=\exp(ih)$. 
For every $t\in\Bbb R$, define the unitary $u_t=\exp(ith)\in L\Bbb Z$. Then Ad$(u_t)\in$Aut$(\tilde B)$, so can define the tensor product automorphism $$\theta_t=\otimes_{\gamma\in\Gamma}(\text{Ad}(u_t))_{\gamma}\in\text{Aut}(\overline{\bigotimes}_{\Gamma}\tilde B),\forall t.$$ Since $\theta_t$ commutes with $\tilde\sigma$, it extends to an automorphism $\theta_t$ of $\tilde M$, for all $t$. Moreover, since $\lim_{t\rightarrow 0}||u_t-1||_2=0$, it follows that $\theta_t\rightarrow$id, as $t\rightarrow 0$, in the pointwise $||.||_2$-topology.

\vskip 0.1in

Note that our proof of Theorem 2
parallels the proof of Lemma 5.2. in [P4]. For this reason we
need to re-establish two main ingredients used there: the
"transversality" of the deformation $\theta_t$ (Lemma 2.1. in [P4])
and the "spectral gap" property  of the inclusion $M\subset\tilde M$ provided by Lemma 5.1. in [P4].

\vskip 0.05in
Firstly,  let $\beta$ be the
automorphism of $\overline{\bigotimes_{\Gamma}}\tilde {B}$ defined by
$\beta_{|\overline{\bigotimes}_{\Gamma}B}=\text{id}_{|\overline{\bigotimes}_{\Gamma}B}$ and by
$\beta ((u)_{\gamma})=(u^*)_{\gamma}$, for every $\gamma\in \Gamma$.
Since $\beta$ commutes with $\tilde{\sigma}$, it extends to an automorphism of $\tilde{M}$, still denoted $\beta $, 
which satisfies ${\beta}^2=\text{id},\beta_{|M}=\text{id}_{M}$ and $\beta\theta_t\beta=\theta_{-t}$, for all $t\in \Bbb{R}$.
As pointed out by Popa (Lemma 2.1 in [P4]), any deformation $\theta_t$ which posseses a symmetry $\beta$
satisfying these relations, automatically verifies the following "transversality" condition:

\proclaim {Lemma 4 ([P4])} For all $t$ and for all $x\in M$ we have that
 $$||\theta_{2t}(x)-x||_2\leq 2||\theta_t(x)-E_{M}(\theta_t(x))||_2.$$
\endproclaim

Secondly, we show that the statement of Lemma 5.1. in [P4] still holds true in our context. 
Although our proof goes
along the same lines as in [P4], there are a few computational 
differences, which we
address below. Briefly, the difference comes from the fact that in [P4], the algebra $\tilde B$ is the tensor product $B\overline{\otimes}B$, rather than the free product $B*L\Bbb Z$, as in our case. This implies that as a $B-B$ bimodule, $L^2\tilde B$ is isomorphic to ${L^2B}^{\oplus\infty}$, in the context of [P4], and, respectively, to ${L^2B}\oplus(L^2B\overline{\otimes}L^2B)^{\oplus\infty}$, in our context.
In turn, this gives different formulas for the $M-M$ bimodule $L^2\tilde M$, depending on the context.

\proclaim {Lemma 5} Assume that $(B,\tau)$ is amenable and let $P\subset M$ be a von Neumann subalgebra with no amenable direct summand. Then 
$P'\cap {\tilde M}^{\omega}\subset M^{\omega}$.
\endproclaim

{\it Proof.} We first prove that the $M-M$
 Hilbert bimodule $L^2\tilde M\ominus L^2M$ is weakly contained in the
 $M-M$ Hilbert bimodule $(L^2M\overline{\otimes}L^2M)^{\oplus\infty}$.

For this, let $\Cal{B}=\{\xi_k\}_{k\in\Bbb N}\subset B$ be an orthonormal
basis for $L^2B$ with $\xi_0=1$. Then  $\tilde
{\Cal{B}}=\{u^{n_1}\xi_{i_1}u^{n_2}\xi_{i_2}\ldots\xi_{i_k}|
n_1,..,n_k\in \Bbb{Z}, i_1,..,i_k\in \Bbb{N}, k\in\Bbb N \}$ is an orthonormal basis for
$L^2\tilde {B}$ which contains $\Cal B$. Since  $\tilde \Cal B\setminus\Cal B$ is infinite,  we can enumerate $\tilde \Cal B\setminus \Cal B=\{\xi_k\}_{k\in\Bbb Z\setminus\Bbb N}$.

Let $\tilde {I}$ (resp. $I$) be the set of sequences
$i=(i_g)_{g\in\Gamma}$ with $i_g\in \Bbb Z$, for all $g\in\Gamma$
(resp. $i_g\in\Bbb{N}$, for all $g\in\Gamma$) such that $\Delta_i:=\{{g\in\Gamma|i_g\neq 0}\}$ is finite. For every $i\in \tilde I$,
define $\eta_i=\otimes_{g} \xi_{i_g}$. Then
$\Cal{\tilde C}=\{\eta_i|i\in\tilde I\}$ and
$\Cal{C}=\{\eta_i|\in
 I\}$  are
orthonormal bases for $L^2(\overline {\otimes}_{\Gamma} \tilde B)$
and respectively for $L^2(\overline{\otimes}_{\Gamma} B)$.

Now, 
let $J$ be the set of  $i=(i_g)_{g\in\Gamma}\in \tilde I$
having the property that for every $i_g\neq 0$ the element
$\xi_{i_g}$ starts and ends with a nonzero power of $u$. Then
it is clear that $$ L^2(\tilde M)\ominus L^2(M)=\bigoplus_{i\in J}
L^2(M\eta_{i}M)\tag 1$$  as Hilbert $M-M$ bimodules.

Note that the Bernoulli action induces an action
of $\Gamma$ on $\tilde I$ by translation and that for every element $i \in \tilde{I}$ its stabilizing group
$\Gamma_{i}\leq\Gamma$ under this action is finite. Since  $\Gamma_i$ satisfies
$\Gamma_{i}(\Gamma\setminus\Delta_{i} )=\Gamma\setminus\Delta_{i}$
we can consider the following crossed product algebra
$K_{i}:=(\overline{\bigotimes}_{\Gamma\setminus\Delta_{i}}B)
\rtimes\Gamma_{i}$.

Next, we claim that for all $i\in J$ we have that  $$L^2(M\eta_{i}M)\cong L^2(\langle
M,K_{i}\rangle,Tr)\tag 2$$  as Hilbert $M-M$ bimodules via the map
$x\eta_{i}y\rightarrow xe_{K_{i}}y$, where as usual  $\langle
M,K_{i}\rangle$  denotes the basic construction corresponding to the
inclusion $K_i\subset M$ and $Tr$ is the canonical trace on it.
To show this it suffices to verify that
$$<(xu_{\gamma})\eta_{i}(yu_{\gamma'}),\eta_{i}>_{\tau}=<
(xu_{\gamma})e_{K_{i}}(yu_{\gamma'}),e_{K_{i}}>_{Tr},\forall x,y\in
\overline{\bigotimes}_{\Gamma} B,\forall \gamma,\gamma'\in\Gamma\tag 3$$
Note that every element of $\overline{\bigotimes}_{\Gamma}B$ can be approximated in $|| \cdot || _2$ by finite linear
combinations of elements from $\Cal{C}$. Thus, in order  to prove identity (3) we can assume that  $x$ and $y$ are of the form $x=\otimes x_{g},y=\otimes y_{g}$, where $x_g,y_g\in B$, for all $g\in\Gamma$ and $x_g=y_g=1$, for all but finitely many $g\in\Gamma$. 
For such $x$ and $y$ the left side of (3) equals $$<(xu_{\gamma})\eta_{i}(yu_{\gamma'}),\eta_{i}>_{\tau}=\delta_{\gamma\gamma',e}\tau(\eta_i^{*}x\tilde\sigma(\gamma)(\eta_iy))=\delta_{\gamma\gamma',e}\prod_{g\in\Gamma}\tau(\xi_{i_g}^{*}x_g\xi_{i_{\gamma^{-1}g}}y_{\gamma^{-1}g})\tag 4$$ Remark that if $\xi_{i_1},\xi_{i_2}\in\Cal{\tilde B}$ are either equal to 1 or start and end with a non-zero power of $u$, then for all $x,y\in B$ we have that $\tau(\xi_{i_1}^*x\xi_{i_2}y)=0$, unless $i_1=i_2$. Moreover, if $i_1=i_2\not=0$, then  $\tau(\xi_{i_1}^{*}x\xi_{i_2}y)=\tau(x)\tau(y)$ and if $i_1=i_2=0$, then  $\tau(\xi_{i_1}^{*}x\xi_{i_2}y)=\tau(xy).$
 Using this remark and (4), we get that left side of (3) equals $$1_{\Gamma_i}(\gamma)\delta_{\gamma\gamma',e}\prod_{g\in\Delta_i}\tau(x_g)\tau(y_{\gamma^{-1}g})\prod_{g\in\Gamma\setminus\Delta_i}\tau(x_gy_{\gamma^{-1}g})\tag 5$$

Now, note that if  $x=\otimes_{g}x_g$ and $\gamma\in\Gamma$, then $E_{K_i}(xu_{\gamma})=1_{\Gamma_i}(\gamma)\prod_{g\in\Delta_i}\tau(x_g)\otimes_{g\in\Gamma\setminus\Delta_i}x_g$. Using this, it is immediate   
that (5) equals $$\tau(E_{K_i}(xu_{\gamma})E_{K_i}(yu_{\gamma'}))=<(xu_{\gamma})e_{K_{i}}(yu_{\gamma'}),e_{K_{i}}>_{Tr}.$$ 
This proves (3) and consequently (2).

Finally, since $B$ is amenable and $\Delta_i$ is  finite, we get that $K_i$ is an amenable von Neumann algebra, for all $i\in J$. This
further implies that $L^2(\langle M,K_{i}\rangle,Tr)\prec
L^2(M)\overline{\otimes} L^2(M)$  and by combining this fact with (1) and
(2), we get the desired weak containment of Hilbert bimodules.
From this it follows, by the argument in the proof of Lemma 2.2. in [P6], that $P'\cap {\tilde M}^{\omega}\subset M^{\omega}$. 
\hfill$\blacksquare$

\vskip 0.1in

We end this section by  recalling {\bf Popa's  technique of intertwining subalgebras} in a finite von Neumann algebra. For this, let $(M,\tau)$ be a finite von Neumann algebra with a normal, faithful trace $\tau$ and let $Q,B\subset M$ be two von Neumann subalgebras. Let $E_B:M\rightarrow B$ denote the  $\tau$-preserving conditional expectation onto $B$. Then we say that a {\bf corner of $Q$ embeds into $B$} if one of the following equivalent conditions (see section 2 in [P1] for the proof of the equivalence) holds true:
\vskip 0.05in
(a) There exist non-zero projections $q\in Q$ and $p\in B$, a normal $*$-homomorphism $\phi:qQq\rightarrow pBp$ and a non-zero partial isometry $v\in M$ such that $xv=v\phi(x)$, for all $x\in qQq$ and $vv^*\in (qQq)'\cap qMq$ , $v^*v\in (\phi(qQq))'\cap pMp$.
\vskip 0.05in
(b) There exist $\varepsilon>0$ and $a_1,..,a_n\in M$ such that $\sum_{i,j=1}^n||E_{B}(a_iua_j^*)||_2^2\geq\varepsilon,$ for all $u\in\Cal U(Q)$.
\vskip 0.05in
Also note that if $B_1,..,B_k$ are von Neumann 
subalgebras such that 
there exist $\varepsilon>0$ and $b_1,..,b_k\in M$ such that $
\sum_{i=1}^k||E_{B_i}(b_iub_i^*)||_2^2\geq \varepsilon,\forall u\in\Cal U(Q)$,
then a corner of $Q$ embeds into $B_i$, for some $i\in\{1,..,k\}$ 
(see the proof of 4.3. in [IPP]).

\vskip 0.1in
\head \S { 3. Proofs of main results.}\endhead
{\bf Proof of Theorem 2.} By working with amplifications we can assume $p=1$. Also, we can assume that $Q$ is diffuse. 
\vskip 0.02in
We first use a spectral gap argument due to Popa to show that $\theta_t$ converge uniformly to id$_{\tilde M}$ on ($Q)_1$ ([P4]). To this end, let $\varepsilon>0$. Since $P$ has no-amenable direct summand,  by Lemma 5 we have that $P'\cap {\tilde M}^{\omega}\subset M^{\omega}$ . Therefore we can find unitaries $u_1,u_2,..,u_n$ and $\delta>0$ such that if $x\in(\tilde M)_1$ satisfies $||[u_i,x]||_2\leq\delta$, for all $i\in\{1,..,n\}$, then $||x-E_M(x)||_2\leq\varepsilon.$
 Let $t_0>0$ such that for all $t\in [0,t_0]$ and all $i$, we have that $||\theta_t(u_i)-u_i||_2\leq \delta/2$. Then we deduce that for all $x\in (Q)_1$, $$||[u_i,\theta_t(x)]||_2\leq 2||\theta_t(u_i)-u_i||_2+||[\theta_t(u_i),\theta_t(x)]||_2\leq \delta,\forall i\in\{1,..,n\}.$$

Thus, by the above  we get that $||\theta_t(x)-E_M(\theta_t(x))||_2\leq\varepsilon,$ for all $x\in (Q)_1$. Lemma 4 then gives that $||\theta_{2t}(x)-x||_2\leq 2\varepsilon,$ for all $x\in (Q)_1$ and all $t\in [0,t_0]$, hence $\theta_t$ converges uniformly to id$_{\tilde M}$ on the unit ball of $Q$.
\vskip 0.02in
Now, recall that Theorem 3.3. in [I] asserts that if $Q\subset M$ is a relatively rigid von Neumann subalgebra
then either (1) a corner of $Q$ embeds either into  $L\Gamma$ or (2) a corner of $Q$ embeds into $\overline{\bigotimes}_{F}B$, for some finite set $F\subset\Gamma$. Since the proof of 3.3. in [I] only uses the  fact that $\theta_t$ converges uniformly id$_{\tilde M}$ on ($Q)_1$, we deduce that either (1) or (2) above are satisfied in our situation.
 
 \vskip 0.02in
In  case (1), the conclusion follows from the proof of Lemma 5.2. in [P4]. So, to end the proof we  need to argue that case (2) leads to a contradiction. Assume therefore that there exist a finite set $F\subset \Gamma$, projections $p\in  \overline{\bigotimes}_F B$, $q\in Q$, a homomorphism $\phi: qQq\rightarrow p(\overline{\bigotimes}_{F}B)p$ and a non-zero partial isometry $v\in M$ such that $vv^*\leq q$ and $xv=v\phi(x)$, for all $x\in qQq$. 

Then $\phi(qQq)$ is a diffuse von Neumann subalgebra of $p(\overline{\bigotimes}_{F}B)p$. Hence if we denote  $K=FF^{-1}$, then Lemma 1.5. in [I] implies that   $$\phi(qQq)'\cap pMp \subset\sum_{\gamma\in K}(\overline{\bigotimes}_{\Gamma}B)u_{\gamma}.$$  Since $P\subset Q'\cap M$, it follows that $v^*Pv\subset \phi(qQq)'\cap pMp$. Thus, $$v^*Pv\subset  \sum_{\gamma\in K}(\overline{\bigotimes}_{\Gamma}B)u_{\gamma},$$ hence  a corner of $P$ embeds into $\overline{\bigotimes}_{\Gamma}B$. This is however a contradiction, since $B$ is amenable, while $P$ has no amenable direct summand.\hfill$\blacksquare$
\vskip 0.1in
{\bf Proof of Corollary 3.}   To prove (i), suppose by contradiction that $P=Q'\cap pMp$ is non-amenable. Then we can find a non-zero projection $z\in\Cal Z(P)$ such that $Pz$ has no amenable direct summand. Since $[Pz,Qz]=0$, by applying Theorem 2, we deduce that there exists a non-zero partial isometry $v\in M$ such that $v^*v\in (Qz)'\cap zMz$ and $vQv^*\subset L\Gamma$. Since $Q$ is diffuse it contains a sequence of unitaries  $\{u_n\}_{n\geq 0}$ which tends weakly to 0.
Then $||vu_nv^*||_2=||vv^*||_2$, for all $n$. On the other hand, since $u_n\in p(\overline{\bigotimes}_{\Gamma}B)p$, for all $n$  and since $u_n\rightarrow 0$ weakly, one easily gets that $||E_{L\Gamma}(vu_nv^*)||_2\rightarrow 0$, as $n\rightarrow\infty$.  This gives a contradiction, thus $P$ has to be amenable.
\vskip 0.02in
For the proof of (ii), suppose that $N'\cap Q^{\omega}\not=\Bbb C1$. Since $N$ is a non-amenable factor, then a result due to Popa (see Lemma 7 in [O1]) implies that we can find a diffuse von Neumann subalgebra $Q_0$ of $Q$ such that $P=Q_0'\cap N$ is non-amenable. Part (i) of this corollary then leads to a contradiction. 
\hfill$\blacksquare$
\vskip 0.1in
  
In the proof of the next result we will use Feldman-Moore's construction of the von Neumann algebra associated to an equivalence relation, which we now recall ([FM]).
Let $\Cal R$ be a countable, m.p. equivalence relation on a probability space $(X,\mu)$. The {\bf full group} of $\Cal R$ (denoted $[\Cal R]$)  is the group of automorphisms $\phi\in$Aut($X,\mu)$ such that $\phi(x)\sim_{\Cal R}x$, a.e. $x\in X$.
On $\Cal R$  consider the measure $\nu$ given by $ \nu(K)=\int_{X}|K\cap {\Cal R}^{x}|\text{d}\mu(x),$ for every $K\subset \Cal R$, and let $L^2( \Cal R,\nu)$ be the associated Hilbert space. For every $\phi\in [\Cal R]$, define the unitary $u_{\phi}\in\Bbb B(L^2(\Cal R,\nu))$ by $$(u_{\phi}g)(x,y)=g(\phi(x),y),\forall g\in L^2(\Cal R,\nu), (x,y)\in \Cal R.$$  Also, represent $L^{\infty}(X,\mu)$ on $L^2(\Cal R,\nu)$ by $$L_f(g)(x,y)=f(x)g(x,y),\forall f\in L^{\infty}(X,\mu),g\in L^2(\Cal R,\nu), (x,y)\in\Cal R.$$ 
The {\bf von Neumann algebra associated to} $\Cal R$ is then defined as the von Neumann algebra generated by $\{u_{\phi}|\phi\in [\Cal R]\}$ and $L^{\infty}(X,\mu)$ and  is denoted $L(\Cal R)$.  Note that $L^{\infty}(X,\mu)$ is a {\bf Cartan subalgebra} of $L(\Cal R)$, i.e. is regular and maximal abelian and that $L(\Cal R)$ is a factor if and only if $\Cal R$ is {\bf ergodic}. Moreover, every central projection $p$ of $L(\Cal R)$ is of the form $p=1_{Y}$, where $Y\subset X$ is an $\Cal R$-invariant measurable subset. In this case, the Cartan subalgebra inclusion associated to $\Cal R_{|Y}$ is isomorhpic to $(L^{\infty}(X,\mu)p\subset L(\Cal R)p)$.

Also,  recall that $\Cal R$ is called  {\bf strongly ergodic} if whenever
 $\{A_n\}_n\subset X$ is a sequence of measurable sets such that
 $\lim_{n\rightarrow\infty}\mu(\phi(A_n)\Delta A_n)=0$, for all
 $\phi\in[\Cal R]$, then $\lim_{n\rightarrow\infty}\mu(A_n)(1-\mu(A_n))=0$ ([Sc], see also [CW],[JS]).  
  By a result of A. Connes, $\Cal R$ is strongly ergodic if and only if $L(\Cal R)'\cap [L^{\infty}(X,\mu)]^{\omega}=\Bbb C1$, where $\omega$ is a free ultrafilter on $\Bbb N$ ([C]).

\vskip 0.05in
\proclaim 
{Proposition 6} Let $\Cal S$ be a countable, m.p.  equivalence relation on a probability space $(X,\mu)$. Then the following are equivalent:
\vskip 0.03in
(1). $Q'\cap L(\Cal S)$ is amenable, for any diffuse von Neumann subalgebra $Q$ of $L^{\infty}(X,\mu)$.
\vskip 0.03in
(2). For any subequivalence relation $\Cal R\subset S$, there exists a partition $\{X_i\}_{i\geq 0}$ of $X$ with $\Cal R$-invariant sets such that (a) $\Cal R_{|X_0}$ is hyperfinite and (b) $\Cal R_{|X_i}$ is ergodic and non-hyperfinite, $\forall i\geq 1$.
\vskip 0.03in
(3). For any subequivalence relation $\Cal R\subset S$, there exists a partition $\{X_i\}_{i\geq 0}$ of $X$ with $\Cal R$-invariant sets such that (a) $\Cal R_{|X_0}$ is hyperfinite and (b) $\Cal R_{|X_i}$ is strongly ergodic, $\forall i\geq 1$.
\vskip 0.03in
(4). (3) holds true for any quotient $\Cal R'$ of a subequivalence relation $\Cal R$ of $\Cal S$.  
\vskip 0.03in 
(5). For any non-atomic probability space $(Y,\nu)$ and for every m.p.,
 onto map $p:X\rightarrow Y$, the equivalence relation $\Cal T=\{(x,y)\in\Cal S|p(x)=p(y)\}$
is hyperfinite.
\vskip 0.03in
(6).  For any non-atomic probability space $(Y,\nu)$, for every m.p.,
 onto map $p:X\rightarrow Y$ and for every hyperfinite equivalence relation $\Cal V$ on $Y$, the equivalence
 relation $\Cal T=\{(x,y)\in\Cal S|p(x)\sim_{\Cal V}p(y)\}$
is hyperfinite.
\endproclaim
{\it Proof.} We first prove the equivalence of the conditions (1)-(4) and then we show 
that (1), (5) and (6) are equivalent.

(1) $\Longrightarrow$ (2). Assume that (1) holds true and let $\Cal R\subset \Cal S$ be a subequivalence relation. Denote by $\Cal Z$ the center of $L(\Cal R)$ and let $p_0\in \Cal Z$ be the maximal projection such that $L(\Cal R)p_0$ is amenable. 

We  claim that $\Cal Z(1-p_0)$ is completely atomic. If not, then we could find a non-zero projection $q$ of $\Cal Z$ such that $q\leq 1-p_0$ and $\Cal Zq$ is diffuse.
Thus, $Q=\Cal Zq\oplus L^{\infty}(X,\mu)(1-q)$ is a diffuse subalgebra of $L^{\infty}(X,\mu)$, so by (1) its relative commutant $Q'\cap L(\Cal S)$ is amenable. In particular, it follows that $L(\Cal R)q$ is amenable, a contradiction to the maximality of $p_0$.
Altogether, we derive that $\Cal Z(1-p_0)$ is completely atomic, hence we can write $\Cal Z(1-p_0)=\oplus_{i\geq 1}\Bbb Cp_i$, for some projections $p_i\in \Cal Z$. 

Now, for every $i\geq 0$, let $X_i\subset X$ be a $\Cal R$-invariant measurable set such that $p_i=1_{X_i}$. Since $L(\Cal R)p_0$ is amenable, Connes-Feldman-Weiss' theorem ([CFW]) implies that $\Cal R_{|X_0}$ is hyperfinite. Also, since $L(\Cal R)p_i$ is a non-amenable factor, we get that $\Cal R_{|X_i}$ is  ergodic and non-hyperfinite, for all $i\geq 1$. 
\vskip 0.03in
(1) $\Longrightarrow$ (3). In the above context, it suffices to show that $\Cal R_{|X_i}$ is strongly ergodic, for all $i\geq 1$.  Assume by contradiction that $\Cal R_{|X_i}$ is not strongly ergodic, for some $i\geq 1$. Then the induced Cartan subalgebra inclusion satisfies 
$[L(\Cal R)p_i]'\cap [L^{\infty}(X,\mu)p_i]^{\omega}\not=\Bbb C1$. Since $L(\Cal R)p_i$ is a non-amenable II$_1$ factor, then, as in the proof of part (ii) of Corollary 3, we can find a diffuse von Neumann subalgebra $Q$ of $L^{\infty}(X,\mu)p_i$, such that $Q'\cap L(\Cal R)p_i$ is non-amenable. This however, violates (1). 
\vskip 0.03in
(1) $\Longrightarrow$ (4).
Let $\Cal R$ be a subequivalence relation of $\Cal S$ and let $\Cal R'$ be a {\bf quotient} of $\Cal R$.  
 Recall that this means that  there exists a measurable, measure preserving, onto map $p:(X,\mu)\rightarrow (X',\mu')$ and a set $N\subset X$ with $\mu(X\setminus N)=0$ such that for all $x\in N$, $p$ is a bijection between the $\Cal R$-orbit of $x$ and the $\Cal R'$-orbit of $p(x)$. Then as noted in [P3] (see 1.4.3.) the von Neumann algebra embedding $\theta:L^{\infty}(X',\mu')\rightarrow L^{\infty}(X,\mu)$ given by $\theta(f)=f\circ p$, for all $f\in L^{\infty}(X',\mu'),$ extends to an embedding $\theta:L(\Cal R')\rightarrow L(\Cal R)$.   Using this observation, the above proof applies verbatim  to show that $\Cal R'$ verifies (3).
\vskip 0.03in
Since it is obvious that (4) $\Longrightarrow$ (3) $\Longrightarrow$ (2), to complete 
the proof of the equivalence of conditions (1)-(4) we only need to show that (2) $\Longrightarrow$ (1). 
To this end, let $Q\subset L^{\infty}(X,\mu)$ be a diffuse von Neumann subalgebra and denote $P=Q'\cap L(\Cal R)$. Since $L^{\infty}(X,\mu)\subset P$, by a result of H. Dye ([D]), we can find a subequivalence relation $\Cal R$ of $\Cal S$ such that $P=L(\Cal R)$. Now, (2) implies that if $\Cal Z$ denotes the center of $L(\Cal R)$ and if  $p_0\in\Cal Z$ is the maximal  projection such that $L(\Cal R)p_0$ is amenable, then $\Cal Z(1-p_0)$ is completely atomic.
On the other hand, since $Q\subset \Cal Z$, we get that $\Cal Z$ is diffuse. 
Hence, we must have that $p_0=1$, therefore $P$ is amenable.
\vskip 0.03in
(1) $\Longleftrightarrow$ (5). Note that every diffuse von Neumann subalgebra
 $Q$ of $L^{\infty}(X,\mu)$ is of the form $L^{\infty}(Y,\nu)\simeq
\{f\circ p|f\in L^{\infty}(Y,\nu)\}$, where
$(Y,\nu)$ is a non-atomic probability space and $p:X\rightarrow Y$ is a measure 
preserving, onto map. For $Q=L^{\infty}(Y,\mu)$, we claim that
$Q'\cap L(\Cal S)=L(\Cal T)$, where $\Cal T=\{(x,y)\in\Cal S|p(x)=p(y)\}$.
Since $L(\Cal T)$ is amenable if and only if $\Cal T$ is hyperfinite ([CFW]), the claim 
implies the equivalence of (1) and (5).
To prove the claim, let $\Cal T'\subset \Cal S$ be a subequivalence relation
such that $Q'\cap L(\Cal S)=L(\Cal T')$ ([D]). Thus, if $\phi$ is an automorphism
of $(X,\mu)$, then $\phi\in [\Cal T']$ if and only if $[u_{\phi},Q]=0$.
 On the other hand, the fact that $u_{\phi}$ commutes 
with $Q=L^{\infty}(Y,\nu)$ is equivalent to $p(\phi(x))=p(x)$, a.e. $x\in X$.
Since the latter condition is in turn equivalent to $\phi\in [\Cal T]$,
we altogether get that $[\Cal T']=[\Cal T]$, thus $\Cal T'=\Cal T$.
\vskip 0.03in
(5) $\Longleftrightarrow$ (6). Since (6) clearly implies (5), we only need to show the converse. For this, assume that (5) holds, let $p:X\rightarrow Y$ be a m.p., onto map and let $\Cal V=\cup_{n}\Cal V_n$ be a hyperfinite equivalence relation, where $\Cal V_n$ are finite equivalence relations. For every $n$, let $\pi_n:Y\rightarrow Y_n:=Y/\Cal V_n$ be the natural projection. Then $\nu_n={\pi_n}_{*}(\nu)$ is a non-atomic probability measure. Since the map $\pi_n\circ p:X\rightarrow Y_n$ is m.p. and onto, by (5) we deduce that the equivalence relation $\Cal T_n=\{(x,y)\in\Cal S|\pi_n(p(x))=\pi_n(p(y))\}=\{(x,y)\in\Cal S|p(x)\sim_{\Cal V_n} p(y)\}$ is hyperfinite, for all $n$. Finally, since $\Cal T=\{(x,y)\in\Cal S|p(x)\sim_{\Cal V}p(y)\}=\cup_{n}\Cal T_n$, we get that $\Cal T$ is hyperfinite.
\hfill$\blacksquare$

\vskip 0.05in
{\bf Proof of Theorem 1.}  Denote $B=L^{\infty}(X,\mu)$ and  on $\overline{\bigotimes}_{\Gamma}B$ we consider the Bernoulli action of $\Gamma$. Then the inclusion $(L^{\infty}([0,1]^{\Gamma})\subset L(\Cal S))$ is naturally  identified with the inclusion $(\overline{\bigotimes}_{\Gamma}B\subset (\overline{\bigotimes}_{\Gamma}B)\rtimes\Gamma)$. Theorem 1 then follows by combining Corollary 3 (i) and Proposition 6.
\hfill$\blacksquare$
\vskip 0.05in
We end this section by noticing a  stronger version of Theorem 1.

\proclaim {Theorem 7} Let $\Gamma\curvearrowright I$ be an action of countable group $\Gamma$ on a set $I$ such that the stabilizer $\Gamma_i=\{\gamma\in\Gamma|\gamma i=i\}$ is amenable, for every $i\in I$.
Let $(X,\mu)$ be a probability space and let $\Cal S_0\subset X\times X$ be a hyperfinite, measure preserving equivalence relation. On $(X,\mu)^{I}$  consider the equivalence relation $\Cal S$ given by: $(x_i)\sim_{\Cal S} (y_i)$ if there exists $\gamma\in\Gamma$ and $F\subset I$ finite such that $x_{i}=y_{\gamma^{-1}i}$, for all $i\in I\setminus F$ and $x_{i}\sim_{\Cal S_0}y_{\gamma^{-1}i}$, for all $i\in F$.  
\vskip 0.02in
 Then $\Cal S$ verifies the conclusion of Theorem 1.
\endproclaim
{\it Proof.} Let $A\subset M$ be the Cartan subalgebra inclusion associated to $\Cal S$. If we denote $B=L(\Cal S_0)$, then it is easy to see that $M$ can be identified with $(\overline{\bigotimes}_{I}B)\rtimes_{\sigma}\Gamma$, where $\sigma$ is the Bernoulli action induced by the  action $\Gamma\curvearrowright I$. Moreover, under this identification $A\subset\overline{\bigotimes}_{I}B$. Also, note that since $\Cal S_0$ is hyperfinite, $B$ is amenable.

To get the conclusion, by the proof of Proposition 6, it suffices to show that if $Q\subset \overline{\bigotimes}_{I}B$ is a diffuse von Neumann subalgebra, then $P=Q'\cap M$ is amenable. 
To this end, we follow the same the lines as in the proof of Theorem 2, leaving some of the details to the reader. 
We start by observing that the context and results of Section 2 extend here.
 Indeed, denote
 $M=(\overline{\bigotimes}_{I}B)\rtimes \Gamma$ and $\tilde M=(\overline{\bigotimes}_{I}\tilde B)\rtimes\Gamma$, where $\tilde B=B*L\Bbb Z$. Then, since
every stabilizer $\Gamma_i$ is  amenable,  Lemma 5 holds true in this context. 
 Moreover, if $\theta_t:\tilde M\rightarrow\tilde M$ is defined analogously (i.e. ${\theta_t}_{|\overline{\bigotimes}_I B}=\otimes_{i\in I}(\text{Ad}(u_t))_i$ and $\theta_t(u_{\gamma})=u_{\gamma},$ for all $\gamma\in\Gamma$), then  Lemma 4 also holds true.

Next, assume by contradiction that  $P$ is non-amenable and let $z\in P$ be a central projection such that $Pz$ has no amenable direct summand. Using the spectral gap argument in the proof of Theorem 2, we deduce that $\theta_t$ converges uniformly to id$_{\tilde M}$ on $(Qz)_1$. Since $Q\subset\overline{\bigotimes}_{I}B$, the proof of Theorem 3.3. in [I] implies that a corner of $Q$ embeds into $\overline{\bigotimes}_{F}B$, for some  finite set $F\subset I$. 

Thus, we find projections $p\in  \overline{\bigotimes}_F B$, $q\in Q$, a homomorphism $\phi: qQq\rightarrow p(\overline{\bigotimes}_{F}B)p$ and a non-zero partial isometry $v\in M$ such that $vv^*\in (qQq)'\cap qMq$ and $xv=v\phi(x)$, for all $x\in qQq$. 
Since $\phi(qQq)\subset p(\overline{\bigotimes}_{F}B)p$ is a diffuse von Neumann subalgebra,  Lemma 1.5. in [I] implies that 
$$\phi(qQq)'\cap pMp \subset\sum_{\gamma\in K}u_{\gamma}(\overline{\bigotimes}_{\Gamma}B)\tag 1$$ where $K=\{\gamma\in\Gamma|\exists i,j\in F, \gamma i=j\}$. 

 Further, since $P\subset Q'\cap M$, it follows that $v^*Pv\subset \phi(qQq)'\cap pMp$. By using this together with (1) we get that  $$v^*Pv\subset  \sum_{i\in F,\gamma\in L}u_{\gamma}[(\overline{\bigotimes}_{\Gamma}B)\rtimes\Gamma_i]\tag 2$$ where  $L\subset\Gamma$ is a finite set such that $K\subset\cup_{i\in F}L\Gamma_i$.
For every $i\in F$, denote by $M_i=(\overline{\bigotimes}_{\Gamma}B)\rtimes\Gamma_i$, by $e_i:L^2(M)\rightarrow L^2(M_i)$  the orthogonal projection and by $E_{M_i}={e_i}_{|M}:M\rightarrow M_i$ the conditional expectation onto $M_i$.
Then, for every $\gamma\in \Gamma$, $u_{\gamma}e_iu_{\gamma}^*$ is the orthogonal projection onto $L^2(u_{\gamma}M_i)$, thus (2) rewrites as $$v^*xv=(\bigvee_{i\in F,\gamma\in L}u_{\gamma}e_iu_{\gamma}^*)(v^*xv),\forall x\in P\tag 3$$	Since the projections $u_{\gamma}e_iu_{\gamma}^*$ mutually commute, we have that $ \bigvee_{i\in F,\gamma\in L}u_{\gamma}e_iu_{\gamma}^*\leq\sum_{i\in F,\gamma\in L}u_{\gamma}e_iu_{\gamma}^*$. By combining this with (3), we get that $$||v^*xv||_2^2\leq\sum_{i\in F,\gamma\in L}||(u_{\gamma}e_iu_{\gamma}^*)(v^*xv)||_2^2=\sum_{i\in F,\gamma\in L}||E_{M_i}(u_{\gamma}^*v^*xv)||_2^2,\forall x\in  P\tag 4$$

Now, recall that $vv^*\in (qQq)'\cap qMq=Pq$, so we can write $vv^*=p'q$, for some projection $p'\in P$. Then $p'E_P(q)=E_P(p'q)=E_P(vv^*)\not=0$, hence we can find a projection $p''\in p'Pp'$ and a constant $C>0$ such that $p''E_P(q)\geq Cp''$. A simple computation then shows that for all $x\in\Cal U(p''Pp'')$, we have that $$||v^*xv||_2^2=\tau(vv^*x^*vv^*x)=\tau(p'qx^*p'qx)=\tau(qx^*x)=\tag 5$$ $$\tau(E_P(q)x^*x)\geq C\tau(p''x^*x)=C\tau(p'').$$

Putting  (4) and (5) together we deduce that $$\sum_{i\in F,\gamma\in L}||E_{M_i}(u_{\gamma}^*v^*xv)||_2^2\geq C\tau(p''),\forall x\in\Cal U(p''Pp'').$$

This in turn implies, that we can find $i\in F$ such that  a corner
 of $P$ embeds into $(\overline{\bigotimes}_{\Gamma}B)\rtimes\Gamma_i$ (by the beginning of the proof of 4.3. in [IPP]). 
Since $(\overline{\bigotimes}_{\Gamma}B)\rtimes\Gamma_i$ is amenable (both $B$ and $\Gamma_i$ are amenable), while $P$ has no amenable direct summand, we get a contradiction.\hfill$\blacksquare$

\head \S {4. Applications.}\endhead 
 
{\bf 4.1. Solid II$_1$ factors.}
The first examples  of solidity in von Neumann algebras are  due to Ozawa who proved that $L\Gamma$ is solid, for every hyperbolic group $\Gamma$ ([O1]). 
Recall in this respect that a II$_1$ factor $M$ is called {\bf solid} if for any diffuse von Neumann subalgebra $A$ of $M$, the relative commutant $A'\cap M$ is amenable. 

\proclaim {Corollary 8} Let $\Gamma$ be an ICC countable group. Then $L\Gamma$ is solid if and only if $(\overline{\bigotimes}_{\Gamma}R)\rtimes\Gamma$ is solid, where $R$ denotes the hyperfinite II$_1$ factor.\endproclaim
{\it Proof.} Assume that $L\Gamma$ is solid and let $P\subset M:=(\overline{\bigotimes}_{\Gamma}R)\rtimes\Gamma$ be a diffuse von Neumann subalgebra. If the commutant $Q=P'\cap M$ is non-amenable, then we can find a non-zero projection $z\in \Cal Z(Q)$ such that $Qz$ has no amenable direct summand.

Since $[Pz,Qz]=0$ and $\Gamma$ is ICC, we can apply Theorem 2, to deduce that there exists a unitary $u$ in $M$ such that $u(Pz\vee Qz)u^*\subset L\Gamma$.
This, however, contradicts the solidity of $L\Gamma$. Thus $Q$ is amenable, hence $M$ is solid.\hfill$\blacksquare$
\vskip 0.05in
Recently, J. Peterson showed that the existence of a proper cocycle into a multiple of the left regular representation  also implies that $L\Gamma$ is solid, for a countable group $\Gamma$ ([Pe]). When combined with Corollary 9, this gives new examples of solid II$_1$ factors.
\vskip 0.1in

{\bf 4.2. Equivalence relations from percolation theory}. We start by recalling how certain percolations on  graphs naturally induce subequivalence relations of an equivalence relation arising from  a Bernoulli action (see [GL] for a reference). 
Let $\Cal G=(V,E)$ be a transitive graph and let $\Gamma\subset$ Aut$(\Cal G)$ be a subgroup which acts transitively on $V$.
Assume that $\Gamma$ acts freely on $E$ with amenable stabilizers, i.e. $\Gamma_e=\{\gamma\in\Gamma|\gamma e=e\}$ is amenable, for every edge $e\in E$. An example of a such a graph is given by the right Cayley graph $\Cal G$ of a countable, finitely generated group $\Gamma$ with respect to a finite generating set $S\subset \Gamma$. More precisely, $\Cal G$ is the graph with vertex set $V=\Gamma$ and edge set $E=\{(\gamma,\gamma s)|\gamma\in\Gamma,s\in S\}$. 

 On $[0,1]^E$ (endowed with the product Lebesque measure $\nu$), consider the Bernoulli  action of $\Gamma$, given by the action of $\Gamma$ on $E$.  Since the stabilizers $\Gamma_e$ are amenable, we get that the
 induced equivalence relation $\Cal S$ satisfies the conclusion of Theorem 1.

Next, let $\pi:[0,1]^E\rightarrow \{0,1\}^E$ be a $\Gamma$-equivariant Borel map. Then $\pi$ gives rise to a subequivalence relation $\Cal R$ of $\Cal S$.
For this,  identify $\{0,1\}^E$ with the set of subgraphs of $E$ and fix $\rho\in V$.  
We then say two points $x,y\in [0,1]^{E}$ are $\Cal R$-equivalent if and only if there exists $\gamma\in \Gamma$ such that
\vskip 0.01in
(1) $\gamma x=y$ and
\vskip 0.01in
(2) $\gamma^{-1}\rho$ and $\rho$ are in the same connected component of $\pi(x)$ (viewed as a subgraph of $E$).
\vskip 0.03in
Now, since $\pi$ is $\Gamma$-equivariant, the push-forward measure $\pi_{*}\nu$ is a $\Gamma$-invariant probability measure on $\{0,1\}^E$, i.e. a {\bf percolation} on $\Cal G$. 
An interesting question  is when does this percolation have  indistinguishable infinite clusters. Equivalently (by [GL]),  when is the restriction $\Cal R_{|U_{\infty}}$ ergodic, where $U_{\infty}$ is the set of all points which have infinite $\Cal R$-classes.

 Specifically, the answer to this question is conjectured to be true for the {\bf free minimal spanning forest} of $\Cal G$ (denoted FMSF($\Cal G$)) ([LPS]). The FMSF is the percolation/equivalence relation induced by the map $\pi:[0,1]^E\rightarrow \{0,1\}^E$ defined as follows: for every $\omega\in [0,1]^E$ label every edge $e\in E$ with $\omega(e)\in [0,1]$ and then define $\pi(\omega)$ to be the set of edges $e$ which are not  maximal (with respect to the labeling) in any cycle containing them. 
Note that it is known by now that for every finitely generated, non-amenable group $\Gamma$ there exists a Cayley graph $\Cal G$ such that a.e. class of FMSF$(\Cal G)$ is a tree with infinitely many ends ([LPS], [T]). More generally, then same is true for every transitive, unimodular graph $\Cal G$ which satisfies $p_c(\Cal G)<p_u(\Cal G)$.  
By combining these facts with Theorem 1 and the proof of Proposition 11 in [GL], we get the following:
\proclaim {Corollary 9} Let $\Gamma$ be a non-amenable, finitely generated group. Then there exists a Cayley graph $\Cal G$ of  $\Gamma$ such that 
  the  FMSF($\Cal G)$
 admits an invariant, measurable partition $\{X_i\}_{i\geq 1} $ of $[0,1]^E$ for which the restriction  
\vskip 0.03in
 $\text{FMSF}(\Cal G)_{|X_i}$ is a (strongly) ergodic, treeable equivalence relation of (normalized) cost $>$ 1, $\forall i\geq 1$.
\vskip 0.03in
\endproclaim   
More generally, Corollary 9 holds true for any transitive, unimodular graph $\Cal G$ such that $p_c(\Cal G)<p_u(\Cal G)$ and that there exists a group  $\Gamma\subset$ Aut$(\Cal G)$   which acts transitively on $V$ and with amenable stabilizers on $E$. 

Recently, D. Gaboriau and R. Lyons used certain equivalence relations coming 
from percolation to show any non-amenable group $\Gamma$ admits $\Bbb F_2$ as a 
"measurable" subgroup ([GL). Also, they suggested that the free minimal spanning forest might be used to derive their result.
We remark that Corollary 9 together with the proof of Proposition 13 in [GL] shows that this is indeed the case. 
\vskip 0.1in

{\bf 4.3. Strong ergodicity vs. spectral gap.} 
While Theorem 1 proves automatic strong ergodicity of certain non-hyperfinite equivalence relations, it is typically quite difficult to prove strong ergodicity for a given equivalence relation.
In  case the equivalence relation is induced by a group action, one usually deduces strong ergodicity by proving spectral gap of the action.  Recall that a measure preserving, ergodic action $\Gamma\curvearrowright^{\sigma} (X,\mu)$ is said to have {\bf spectral gap}  if the induced representation of $\Gamma$ on $L^2(X,\mu)\ominus \Bbb C1$ has spectral gap.

The  notions of strong ergodicity and spectral gap are however not equivalent in general ([Sc],[HK]). Nevertheless, it is an interesting problem to find classes of group actions for which these notions coincide.
 This is the case if  $\sigma$ is a generalized Bernoulli action ([KT]) or if $\sigma$ comes from an embedding of $\Gamma$ as a dense subgroup of a compact group $G$ ([AN]). 
 Below, we note that the argument in [AN] (Lemma 6) only uses the fact that $\sigma$ has large commutant, thus rendering the following: 
 
\proclaim {Lemma 10 } Let $\Gamma\curvearrowright^{\sigma} (X,\mu)$ be a  m.p. action. Assume that the commutant of $\Gamma$ in {\text Aut}$(X,\mu)$ acts ergodically on $(X,\mu)$. Then $\sigma$ has spectral gap if and only if $\sigma$ is strongly ergodic.
\endproclaim
{\it Proof.} Denote by $\Lambda$ the commutant of $\Gamma$ in Aut$(X,\mu)$.
 We claim that for every measurable sets $A,B\subset X$ and every $\varepsilon> 0$ we can find $\theta\in\Lambda$ such that $ \mu(\theta(A)\cap B)\leq (1+\varepsilon)\mu(A)\mu(B)$. This claim is folklore, but we include a proof for the sake of completeness.
 Assuming that the claim is false, then $$\mu(\theta(A)\cap B)>(1+\varepsilon) \mu(A)\mu(B),\forall\theta\in\Lambda.$$ Thus, if $K$ denotes the 
$||\cdot||_2$-closure of the convex hull of the set $\{1_{\theta(A)}|\theta\in\Lambda\}\subset L^2(X,\mu)$, then $\int_{B}f\text{d}\mu>\mu(A)\mu(B)$ and $\int_{X}f\text{d}\mu=\mu(A)$, for all $f\in K$. Let $f\in K$ be the element of minimal $||\cdot||_2$. Then, since $K\ni g\rightarrow \theta(g)\in K$ is a $||\cdot||_2$-preserving map, for all $\theta\in\Lambda$, and since $f$ is unique, we get that $f$ is $\Lambda$-invariant.
Using the fact that $\Lambda$ acts ergodically on $(X,\mu)$ we get that $f\in \Bbb C1$. Thus $f=\mu(A)1_X$, which contradicts the inequality $\int_{B}f\text{ d}\mu>(1+\varepsilon)\mu(A)\mu(B)$. 						

Equivalently, the claim shows that for every $A,B\subset X$ and every $\varepsilon >0$ there exists $\theta\in\Lambda$ such that $\mu(\theta(A)\cup B)\geq (1-\varepsilon)-(1-\mu(A))(1-\mu(B))$. Let $A \subset X $ be a measurable set. By using the claim  we can  inductively find (as in [AN]) a sequence $\theta_1=1,\theta_2,..,\theta_k...\in\Lambda$ such that $\mu(\cup_{i=1}^{k}\theta_i(A))\geq 1-(1-\mu(A))^{k-1},$ for all $k$. 

Now, assume that  $\sigma$ is strongly ergodic but does not have spectral gap. Thus, we can find measurable sets $\{A_n\}\subset X$ such that $\mu(A_n)>0$, $\lim_{n\rightarrow\infty}\mu(A_n)=0$ and  $\lim_{n\rightarrow\infty}\mu(\gamma A_n\Delta A_n)/\mu(A_n)=0,$ for all $\gamma\in\Gamma$. 
For every $n$, define $k_n=[1/(2\mu(A_n))]+1$. Then by the above, we can find  $\theta_{n,1},..,\theta_{n,k_n}\in\Lambda$ such that $B_n:=\cup_{i=1}^{k_n}\theta_i(A_n)$ verifies $\mu(B_n)>1-(1-\mu(A_n))^{k_n-1}$. 
Then, since $\mu(A_n)\rightarrow 0$, it is easy to check that for a large enough $n$ we  have that $$1-1/\sqrt{e}\leq\mu(B_n)\leq 2/3\tag 1$$
Also, using the fact that $\Gamma$ and $\Lambda$ commute,  it follows that  $$\mu(\gamma B_n\Delta B_n)\leq\sum_{i=1}^{k_n}\mu(\gamma(\theta_i(A_n))\Delta\theta_i(A_n))=k_n\mu(\gamma A_n\Delta A_n)\leq\tag 2$$ $$ \mu(\gamma A_n\Delta A_n)/(2\mu(A_n)),\forall\gamma\in\Gamma.$$

Finally, (1) and (2) contradict the assumption that $\sigma$ is  strongly ergodic.\hfill$\blacksquare$

\vskip 0.1in
Remark that Lemma 11 implies that if the commutant of $\sigma$ acts weakly mixing on $(X,\mu)$ then $\sigma$ has double spectral gap if and only if $\sigma$ has double strongly ergodic (see [P4] for definitions).

\head  References\endhead
\item {[AN]} M. Ab\'ert, N.Nikolov: {\it The rank gradient from a combinatorial point of  view}, preprint math.GR/0701925.  
\item {[C]} A. Connes: {\it Outer conjugacy classes of automorphisms of factors}, Ann. Ec. Norm. Sup. {\bf 8} (1975), 383-419. 

\item {[CFW]} A. Connes, J. Feldman, B. Weiss: {\it An amenable equivalence relations is generated by a single transformation}, Ergodic Th. Dynam. Sys. {\bf 1}(1981), 431-450.

\item {[CW]} A. Connes, B.Weiss: {\it Property (T) and asymptotically invariant sequences}, Israel J. Math. {\bf 37}(1980), 209-210.
\item {[D]} H. Dye: {\it On groups of measure preserving transformations} II, Amer. J. Math. {\bf 85}(1963), 551-576.
\item {[FM]} J. Feldman, C.C. Moore: {\it Ergodic equivalence relations, cohomology, and von Neumann algebras, II}, Trans. Am. Math. Soc. {\bf 234}(1977), 325-359.
\item {[GL]} D. Gaboriau, R. Lyons: {\it A-Measurable-Group-Theoretic Solution to von Neumann's Problem}, math.GR/07111643.
\item {[HK]} G. Hjorth, A. Kechris: {\it Rigidity theorems for actions of product groups and countable Borel equivalence relations}, Memoirs of the Amer. Math. Soc., {\bf 177}, No. 833, 2005.
\item {[I]} A. Ioana: {\it Rigidity results for wreath product II$_1$ factors}, Journal of Functional Analysis {\bf 252}(2007) 763-791.
\item {[IPP]} A. Ioana, J. Peterson, S. Popa: {\it Amalgamated Free Products of w-Rigid Factors and Calculation of their Symmetry Groups}, math.OA/0505589, to appear in Acta Mathematica. 
    \item {[JS]} V.F.R. Jones, K. Schmidt: {\it Asymptotically invariant sequences and approximate finiteness,}
  Amer. J. Math.  {\bf 109}(1987), 91--114.   
\item {[KT]} A. Kechris, T. Tsankov: {\it Amenable actions and almost invariant sets}, preprint 2006, to appear in the Proceedings of the AMS.
\item {[LPS]} R. Lyons, Y. Peres, O. Schramm: {\it Minimal spanning forests,} Ann. Probab., {\bf 27} (2006), 1665-1692.
\item{[MvN]} F. Murray, J. von Neumann: {\it Rings of operators}, IV, Ann. Math. {\bf 44}(1943), 716-808.
\item {[O1]} N. Ozawa: {\it Solid von Neumann algebras}, Acta Math., {\bf 192} (2004), 111-117.
\item {[O2]} N. Ozawa: {\it A Kurosh type theorem for type II$_1$ factors,} Int. Math. Res. Not. (2006) Vol. {\bf 2006}, article ID 97560.
\item {[Pe]} J. Peterson: {\it L$^2$-rigidity in von Neumann algebras}, preprint math.OA/0605033.
\item {[P1]} S. Popa: {\it Strong rigidity of II$_1$ factors arising from malleable actions of w-rigid groups
 I}, Invent. Math. {\bf 165}(2006), 369--408. 
\item {[P2} S. Popa: {\it Strong rigidity of II$_1$ factors arising from malleable actions of
 w-rigid groups II}, Invent. Math. {\bf 165}(2006), 409--451. 
\item {[P3]} S. Popa: {\it Cocycle and orbit equivalence superrigidity for Bernoulli 
actions of Kazhdan groups}, Invent. Math, {\bf 170}(2007), 243-295.
\item {[P4]} S. Popa: {\it On the superrigidity of malleable actions with spectral gap},  Journal of the AMS 
(math.GR/0608429).
\item {[P5]} S. Popa: {\it Deformation and rigidity for group 
actions and von Neumann algebras},  International Congress of Mathematicians. 
Vol. I,  445--477, 
Eur. Math. Soc., Zürich, 2007.
\item {[P6]} S. Popa: {\it On Ozawa's property for free group factors}, Int. Math. Res. Not. (2007) Vol. {\bf 2007}, article ID rnm036. 
\item {[Sc]} K. Schmidt: {\it Amenability, Kazhdan's property T, strong ergodicity and invariant means for ergodic group-actions}, Ergod. Th. Dynam. Sys. {\bf 1} (1981), 223-236.
\item {[T]} \'A. Tim\'ar: {\it Ends in minimal spanning  forests}, Ann. Probab., {\bf 34} (2006), 865-869.

\enddocument